\documentclass[12pt]{amsart}
\usepackage{amsmath}

\newtheorem{thm}{Theorem}
\newtheorem{lem}[thm]{Lemma}
\newtheorem{prop}[thm]{Proposition}

\theoremstyle{definition}
\newtheorem{defn}{Definition}[section]

\theoremstyle{remark}

\theoremstyle{plain}

\newcommand{\Z}{{\mathbf Z}}

\newcommand{\R}{{\mathbf R}}

\newcommand{\supp}{\operatorname{supp}}

\newcommand{\CN}{{\mathcal N}}

\newcommand{\tq}{\tilde q}
\newcommand{\SF}{{\mathcal F}}
\newcommand{\meas}{\mbox{meas}}
\newcommand{\rad}{\operatorname{rad}} 

\numberwithin{equation}{section}

\begin{document}

\title[Spacings between fractional parts of $n^2\alpha$] 
{The distribution of spacings between the fractional parts of $n^2 \alpha$}
\author{Ze\'ev Rudnick$^{(1)}$, Peter Sarnak$^{(2)}$ and Alexandru Zaharescu}

\address{Raymond and Beverly Sackler School of Mathematical Sciences,
Tel Aviv University, Tel Aviv 69978, Israel 
({\tt rudnick@math.tau.ac.il})}

\address{Department of Mathematics, Princeton University, Fine Hall,
Washington Road, Princeton, NJ, 08544, USA
({\tt sarnak@math.princeton.edu})}

\address{Department of Mathematics, University of Illinois, 
1409 West Green Street, Urbana, IL 61801, USA  
({\tt zaharesc@math.uiuc.edu})}

\date{\today} 

\thanks{(1) Supported in part by a grant from 
the U.S.-Israel bi-national Science Foundation} 
\thanks{(2) Supported in part by a grant from 
the U.S.-Israel bi-national Science Foundation and the NSF}


\maketitle

\section{Introduction}

Fix an irrational number $\alpha$. The problem of the distribution
of the local spacings between the members of the sequence 
$n^2\alpha \mod 1$,
$1 \le n \le N,$
has received attention recently (see \cite{B-T,Ca,R-S}).
It arises for example in the study of the local spacing distributions 
between the eigenvalues  of special Hamiltonians.
We order the above numbers
in $[0,1)$ as
\begin{equation}\label{i1}
0 \le \beta _1 \le \beta _2 \le \cdots \le \beta_N < 1
\end{equation}
and set $\beta_{N+j}=\beta_j$.
The $k$-th consecutive spacing measure is defined to be the 
probability measure on $[0,\infty  )$ given by
\begin{equation}\label{i2}
\mu_k (N,\alpha ):=\frac 1N\sum_{j=1}^N \delta_{N(\beta_{j+k}-\beta_j)}
\end{equation}
where $\delta _x$ is a unit delta mass  at $x$.
The problem is to understand the behavior of these measures as 
$N \rightarrow \infty$      
and in particular their dependence on the diophantine approximations to 
$\alpha$.

We say $\alpha $ is of type $K$ if there is $c_{\alpha }>0$ such that
$|\alpha -\frac aq|\ge c_{\alpha }q^{-K }$
for all relatively prime   integers $a$ and $q$. It is easy to see that
if $\alpha $ is not of type 3 then there is a subsequence $N_j
\rightarrow \infty$ for which the measures $ \mu_k (N_j,\alpha )$            
converge to a measure supported on ${\bf N}=\{0,1,2,\cdots\}$. On the
other hand numerical
experiments  \cite{Ca} indicate that for $ \alpha = \sqrt 2$
these $k$-th consecutive spacings behave like what one typically
gets for spacings when placing $N$ numbers in $[0,1]$ uniformly
and independently at random  \cite{Fe}. That is 
$ \mu_k (N,\sqrt 2)$ appears to converge to 
$\mu_k :=\frac {x^k}{k!}e^{-x}dx$.

A standard approach to the analysis of the consecutive spacing
 is via local $m$-level 
correlations.
These are defined as follows: As test functions we use functions
$f(x_1,x_2,\cdots ,x_m)$
which are symmetric in $(x_1,x_2,\cdots ,x_m)$
and which are functions of the difference of the coordinates,
that is $f(x+(t,t,\cdots ,t))=f(x)$ for all $t \in {\bf R}$.
We assume further that $f$ is local, that is 
it is compactly supported modulo the
diagonal. We will call these  {\em admissible} test functions. 
Define the correlations
\begin{equation}\label{i3}
R^{(m)}(N,\alpha,f):=\frac 1N \sum_{1 \le j_1<...<j_m \le N}
f(N(\beta_{j_1},\cdots ,\beta_{j_m})).
\end{equation}
Note that $R^{(m)}$  is not a {\it probability} density and it may 
well tend to infinity as $N \rightarrow \infty$
(think of the case when $\alpha$ is rational).
In  the case that the $\beta$'s in \eqref{i1} come from a random choice of
points in $[0,1)$ these correlations satisfy 
\begin{equation}\label{i4}
R^{(m)}(N,f) \rightarrow \int_{0 \le x_2 \le ...\le x_m
}f(0,x_2,\cdots,x_m)dx_2\cdots dx_m.
\end{equation}

We say that $n^2\alpha \mod 1$, $n \le N$, is 
Poissonian if for all $m \ge 2$ and $f$ as above
\begin{equation}\label{i5}
R^{(m)}(N,\alpha,f) \rightarrow \int_{0 \le x_2 \le ...\le x_m
}f(0,x_2,\cdots,x_m)dx_2\cdots dx_m
\end{equation}
as $N \rightarrow \infty$. 

As with the method of moments in convergence of measures, if the 
$m$-level correlations are Poissonian
then the consecutive spacing measures $\mu_k(N,\alpha)$ converge
to $\mu_k$. Thus Poissonian in the sense of (\ref{i5}) (i.e. for
correlations) implies that as far as local spacings go, the numbers
behave randomly. We will also consider cases where (\ref{i5})
holds along a subsequence $N_j \rightarrow \infty$, in such a 
case we say that 
$n^2\alpha \mod 1$, $1 \le n \le N$ is 
Poissonian on a subsequence.

The results below lead us to the following  

\vskip .5cm
\noindent{\bf Conjecture:}
{\em If $\alpha$ is of type $2+\epsilon$ for every $\epsilon > 0$
and the convergents $\frac aq$ to $\alpha$ satisfy 
$\lim_{q \rightarrow \infty} \frac {\log \tilde q}{\log q} =1$,
where $\tilde q$ is the square free part of $q$, then 
 $n^2\alpha \mod 1$ is Poissonian.}
\vskip .5cm

We note that almost all $\alpha$ (for Lebesgue measure) satisfy the
hypothesis in the Conjecture and that
assuming some standard conjectures in diophantine analysis any
real algebraic irrationality  satisfies these 
hypotheses (see Appendix~\ref{app:A}).

Unfortunately the methods of this paper appear not to be powerful enough
to prove anything for the numbers $\alpha$  in the Conjecture.
They require that $\alpha$ have somewhat better approximations 
by rationals. One of our main results gives conditions on the
diophantine approximations to $\alpha$ which ensure that
 $n^2\alpha \mod 1$ is Poissonian along a subsequence. 
In particular this allows us to conclude that for the topologically
generic $\alpha$ (i.e. in the sense of Baire)  $n^2\alpha \mod 1$ 
is Poissonian along a subsequence. On the other hand the naive
expectation that for any irrational $\alpha$, 
 $n^2\alpha \mod 1$ is Poissonian along a subsequence, fails
dramatically. The source of this phenomenon is large square
factors in the denominator of the convergents to $\alpha$.
We will exhibit an $\alpha$  for which the $5$-level correlations go
to infinity as $N\to \infty$. We also provide an $\alpha$
of type less than three and a sequence of integers
$\{N_j\}_{j=1}^\infty$ along which   
the 5-level correlations diverge to infinity. 

The precise statements are as follows: 
\begin{thm} \label{thm:1}
Let $\alpha \in {\bf R}$. Suppose there are infinitely many rationals 
${b_j}/{q_j}$, with $q_j$ {\em prime},  satisfying
$$
\Big |\alpha - \frac {b_j}{q_j}\Big | < \frac 1{q_j^3} \;.
$$
Then there is a subsequence $N_j \rightarrow \infty $      with
$\frac {\log N_j}{\log q_j} \rightarrow 1$
for which (\ref{i5}) holds for all 
$m \ge 2$ and all $f$. That is to say 
 $n^2\alpha \mod 1$ is Poissonian along this subsequence.
\end{thm}

With a lot more work concerning the exponential sums 
discussed in Section~\ref{sec:3}, for general moduli $q$, we can relax the
condition that $q_j$ be prime in Theorem~\ref{thm:1}. In fact we can prove
the following (we do not go into the proof in this paper)
which shows that for such approximants the size of the square 
free parts $\tilde q_j$ of $q_j$ is decisive.

\vskip .5cm

\noindent {\bf Theorem 1'}: 
{\em Let $\alpha$ be an irrational for which there are infinitely many
rationals  ${b_j}/{q_j}$ satisfying}
$$
\Big |\alpha - \frac {b_j}{q_j}\Big | < \frac 1{q_j^3} \;.
$$
{\em Then the following are equivalent :

(i) There is a subsequence $N_j \rightarrow \infty $      with
$\frac {\log N_j}{\log q_j} \rightarrow 1$
such that 
 $n^2\alpha \mod 1$ is Poissonian along $N_j$.

(ii) }
$$
\lim_{j \rightarrow \infty} \frac {\log \tilde q_j}{\log q_j} =1 \;.
$$

\vskip .25cm

As to the divergence of correlations we have:
\begin{thm}\label{thm:2}
(a) There is an irrational $\alpha$ and a test function  
$f$ such that
$$
R^{(5)}(N,\alpha,f) \gg N^{\frac {1}{4}}(\log N)^{-5},\text{ as } \quad
N \rightarrow \infty \;. 
$$
(b) For every $\sigma>23/8$ there is an $\alpha$  of type $\sigma$ 
and a test function  $f$ such that 
$$
\overline {\lim_{N \rightarrow \infty}} R^{(5)}(N,\alpha,f) = \infty \;.
$$
\end{thm}

The test functions $f$  in Theorem~\ref{thm:2} are nonnegative 
and are supported in a neighborhood of $0$ (modulo the 
diagonal) and the 
source of the divergence is that there are zero density,
 but non-negligible,
clusters among the numbers $n^2\alpha \mod 1$ , $n \le N$.    

We note that these clusters which spoil the correlations
do not have the same effect on the probability measures 
$\mu_k (N,\alpha)$. So it is quite possible for example that the
$\alpha$ in part (b) of Theorem~\ref{thm:2} has its $\mu_k (N,\alpha)$ 
measures converge to the Poissonian $\mu_k$. We have chosen in this
paper to call  $n^2\alpha \mod 1$  Poissonian if the strongest
behavior holds - that is the correlations are Poissonian.

The proofs of the above theorems are based on 
the following closely related
diophantine problem: Consider the spacing distributions (normalized 
to mean spacing $1$ as 
before) of the numbers $\{bn^2/q\}$, $n\leq N$, or what is the same, the
spacing distribution of the integers 
\begin{equation}\label{i6}
n^2b \mod q , \qquad 1 \le n \le N.
\end{equation}
Here $q$ is prime ($q  \rightarrow \infty$),
$b$ is any number not divisible by $q$
and  $N$ is in the range $[q^{1/2+\epsilon} , \frac q{\log q}]$
for some 
$\epsilon  > 0$.
The reason for this range for $N$ is that if $N \le \sqrt q$ and say
$b=1$ 
then the spacing distributions
may be easily determined (since $n^2 < q$ for $n \le N \le \sqrt q$)
and are certainly nonrandom.
 Similarly if $N=q$  then the sequence in (\ref{i6}) consists of 
all the 
quadratic residues (or non-residues) and hence  the spacings are 
integers and so cannot follow a Poissonian law.
In fact the limiting
spacing distributions of $\mu_k(q,q,b)$  were determined by Davenport
\cite{Davenport-quadratic, Davenport-powers}.  
So it is only in the range $N \in [q^{1/2+\epsilon} , \frac q{\log q}]$     
that we can hope for randomization.
The following Theorem shows that indeed, to a certain extent, this
is the case.

Let $R^{(m)}(N,b/q,f)$
denote the scaled $m$-level correlations for the sequence (\ref{i6}).

\begin{thm}\label{thm:3}
Fix  $m \ge 2$ ,  $f$ and $\delta > 0$. Then as 
$q \rightarrow \infty$, $q$ prime
$$
R^{(m)}(N,b/q,f)\rightarrow \int_{0 \le x_2 \le ... \le x_m} 
f(0,x_2,\cdots,x_m) dx_2...dx_m
$$
uniformly for $(b,q)=1$ and
$q^{1-\frac 1{2m}+\delta} \le N \le \frac q{\log q} $.
\end{thm} 

A crucial ingredient in our proof of Theorem~\ref{thm:3} is the Riemann
Hypothesis for curves (of arbitrary large 
genus) over finite fields (Weil \cite{We}).

In the range in which Theorem~\ref{thm:3} applies it gives Poisson 
statistics and Theorem~\ref{thm:3} easily yields Theorem~\ref{thm:1}.
For $m \ge 3$ it is not possible to extend the range of $N$
in Theorem~\ref{thm:3} much further. The reason is related to the previous
divergence  of correlations phenomenon.
For suitable $b$ (depending on $q$) there will be large clusters among
the numbers  $n^{2} b \pmod q$ , $n \le N$.
This is highlighted by the following Theorem.

\begin{thm}\label{thm:4} 
Fix $m \ge 3$ and $\delta > 0$. Then there is a test function $f$ 
such that for   $q^{\frac 12} \le N \le q^{\frac m{m+2}-\delta }$, 
$$
\lim_{q \rightarrow \infty} \max _{(b,q)=1} R^{(m)}(N,b/q,f) = \infty \;.
$$
\end{thm}

{\bf Acknowledgment:} We would like to thank E. Bombieri for his help with 
the application of the  ABC conjecture described in Appendix~\ref{app:A}.
We also thank the referee for  suggesting a stronger version of 
Theorem~\ref{thm:2}~(a) with a simpler proof than our original one.

\section{A comparison lemma} 

We will need to deal with the following situation: We are given two
families of sequences $\CN = \{x_{N}(n):n\leq N\}$ and 
$\CN'=\{x'_{N}(n):N\leq N\}$ in
$[0,1)$ and we wish to compare the 
limiting correlation functions of these two families, 
seeking to show that if 
the correlations exist for one sequence then they exist for the other,
or they diverge for one if they do for the other. 
We show that it can be done if the two sequences are close in a
suitable sense.   
We define the scaled distance between the sequences to be 
$$
\epsilon(\CN, \CN'):= N \max_{n\leq N} |x_{N}(n)- x'_{N}(n)| \;.
$$
A general method for carrying out the comparison is formalized in the
following: 
\begin{lem}[Comparison Lemma]\label{comparison:lem} 
Assume that $\CN,\CN'\subset [0,1)$ are two families of sequences with
$\epsilon(\CN, \CN') \to 0$ as $N\to \infty$. Then for all smooth test 
functions $f$, we have 
$$
\left| R^{(k)}(\CN,f) - R^{(k)}( \CN',f)\right| \leq  
R^{(k)}(\CN,f_+) \epsilon(\CN, \CN')
$$
for $N$ sufficiently large, 
where $f_+\geq 0$ is a smooth admissible test function (depending only
on $f$). 
\end{lem}
\begin{proof}
For notational simplicity, we will do the case of pair correlation
($k=2$). Our test function $f\neq 0$ can then be written as
$f(x_1,x_2)=g(x_1-x_2)$ for some $g\in C_c^\infty(\R)$, say $g$
supported inside $[-\rho,\rho]$. Let $g_+\geq 0$ be smooth, compactly
supported and such that $g_+$ is constant  on $[-2\rho,2\rho]$, where
it equals  $\max|g'|$.  
Set  $f_+(x_1,x_2):=2 g_+(x_1-x_2)$. 
For further notational simplicity also set 
$\delta_{m,n}:=x_N(m)-x_N(n)$ and
$\delta'_{m,n}:=x'_N(m)-x'_N(n)$. 

By the mean value theorem we have 
\begin{equation*}
\begin{split}
R^{(2)}(\CN,f) - R^{(2)}(\CN',f) &= \frac 1N \sum_{1\leq m< n\leq N} 
g(N\delta_{m,n}) - g(N\delta'_{m,n})\\
&=
\frac 1N \sum_{1\leq m< n\leq N} 
g'(N\xi_{m,n}) \cdot N (\delta_{m,n}-\delta'_{m,n})
\end{split}
\end{equation*}
where $\xi_{m,n}$ lies between $\delta_{m,n}$ and
$\delta'_{m,n}$. 

For the difference $R^{(2)}(f,\CN) - R^{(2)}(f,\CN')$
to contain a nonzero contribution from the  term indexed by the pair
$(m,n)$, we must
have at least one of $N\delta_{m,n}$ or $N\delta'_{m,n}$ lying in
$\supp g\subset [-\rho,\rho]$. 
Now $N\xi_{m,n}$ is within $2\epsilon(\CN,\CN')$ of both $N\delta_{m,n}$
and $N\delta'_{m,n}$, which implies that {\em both} lie in
$[-2\rho,2\rho]$, as does $\xi_{m,n}$ if $N$ is sufficiently large so
that $2\epsilon(\CN,\CN')<\rho$. 
Since $g_+$ is constant on $[-2\rho,2\rho]$ we find that 
$g_+(\xi_{m,n}) = g_+(N\delta_{m,n})$. 

Thus we get 
\begin{equation*}
\begin{split}
| R^{(2)}(\CN,f) - R^{(2)}(\CN',f)| &\leq \frac 1N \sum_{1\leq m< n\leq N} 
 g_+(N\delta_{m,n}) \cdot 2\epsilon(\CN,\CN')\\ 
& = R^{(2)}(\CN,f_+) \epsilon(\CN,\CN')
\end{split}
\end{equation*}
as required. 
\end{proof}

\section{Derivation of Theorem~\ref{thm:1}}\label{sec:2}

As an immediate application of the comparison lemma, we derive
Theorem~\ref{thm:1} from Theorem~\ref{thm:3}. 

Fix  $\alpha$. Suppose there are infinitely many rationals 
$ {b_j}/{q_j}$ with $q_j$ prime,  satisfying 
\begin{equation*}\label{e21}
|\alpha -\frac {b_j}{q_j}|<\frac {1}{q_j^3}.
\end{equation*}
We let $N_j=[\frac {q_j}{\log q_j}]$, where [.] denotes the integer
part function. 
Fix an $m\ge 2$ and a test function $f$ as above. We need to show that
\begin{equation*}\label{e22}
\lim_{j\rightarrow \infty }R^{(m)}(N_j,\alpha ,f)=\int _{0\le x_2\le \cdots
\le x_m}f(0,x_2,\cdots ,x_m)dx_2\cdots dx_m.
\end{equation*}
By Theorem~\ref{thm:3} applied to $q=q_j, b=b_j$ and $N=N_j$ we know that 
\begin{equation*}\label{e23}
\lim_{j\rightarrow \infty }R^{(m)}(N_j,b_j/q_j ,f)=\int _{0\le x_2\le \cdots
\le x_m}f(0,x_2,\cdots ,x_m)dx_2\cdots dx_m
\end{equation*}
We use the comparison principle to estimate the difference:
$$
|R^{(m)}(N_j,\alpha ,f)-R^{(m)}(N_j,b_j/q_j,f)| \;. 
$$
Take $\CN'_j =\{\{\alpha n^2\}: n\leq N_j\}$ and $\CN_j =
\{\{b_jn^2/q_j\}: n\leq N_j\}$. 
By lemma~\ref{comparison:lem}, 
\begin{equation}\label{Apply comparison 1} 
|R^{(m)}(\CN'_j,f)-R^{(m)}(\CN_j,f)|  \leq
 R^{(m)}(\CN_j,f_+) \epsilon(\CN_j,\CN_j')
\end{equation}
 for some  admissible test function $f_+ \geq 0$. 
We have 
$$
| \{\alpha n^2\}- \{\frac{b_j}{q_j} n^2 \} | \leq |\alpha -\frac {b_j}{q_j}| n^2
\leq \frac{N_j^2}{q_j^3} \sim \frac 1{N_j (\log N_j)^3}
$$
and thus 
$$
\epsilon(\CN_j,\CN_j') = N_j \max_{n\leq N_j}| \{\alpha n^2\}-
\{\frac{b_j}{q_j} n^2 \}| 
\leq \frac 1{(\log N_j)^3} \to 0 \;.
$$
By Theorem~\ref{thm:3},  $R^{(m)}(\CN_j,f_+)$ is bounded (it converges
as $j\to\infty$).  
Thus we use \eqref{Apply comparison 1} to deduce that 
$$
|R^{(m)}(N_j,\alpha ,f)-R^{(m)}(N_j,b_j/q_j,f)| \to 0
$$
which gives Theorem~\ref{thm:1}. \qed

\section{A divergence principle}\label{sec:4}
We present a mechanism that ensure divergence of high correlations of
the sequence $\{bn^2/q\}$:   The presence of larges square factors in $q$. 

\begin{lem}\label{blowup:lem}
Let $q=uv^2$  with $v>q^{\delta } $ for some 
$\delta >0$, let $\eta>1-\delta$ and suppose that $\log N/\log q >\eta$. 
Let $f\geq 0 $ be a  {\em positive} admissible test
function which is non-vanishing at the origin. 
Then for all $b$, 
\begin{equation}\label{blowup:eq}
 R^{(m)}(N,\frac bq, f) \gg \frac 1N f(0) (\frac {Nv}q)^m \;. 
\end{equation}
In particular 
$R^{(m)}(N,b/q, f)$ will diverge to infinity
for $m$ sufficiently large 
in terms of $\delta $ and $\eta$.
\end{lem}
\begin{proof}
Write $f(x_1,\dots ,x_m) = g(x_1-x_2,\dots x_{m-1}- x_m)$ for $g\in
C_c(\R^{m-1})$, $g\geq 0$, $g(0)\neq 0$. Then 
$$ 
R^{(m)}(N,\frac bq, f) = \frac 1N \sum_{1\leq n_1<\dots <n_m\leq N} 
g(\dots, N\{\frac {bn_j^2}{q}\} -N\{ \frac {b n_{j+1}^2}{q}\}, \dots)
$$
Since $g\geq 0$, we may count only the contribution of those 
$(n_1,\dots ,n_m) $ ($n_j$ distinct) for which all the 
components $n_1,\dots ,n_m$ are divisible by $uv$. 
There are $\gg [N/uv]^m = [Nv/q]^m$ such $m$-tuples.  
If $n=uvn'$ then 
since $q=uv^2$ we have 
$$
\{ \frac {b n^2}{q}\} =\{ bu (n')^2 \} = 0
$$
and so we find 
$$
 R^{(m)}(N,\frac bq, f) \gg \frac 1N f(0) (\frac {Nv}q)^m \;. 
$$
Since $v>q^\delta$ and $N\gg q^\eta$ with $\eta>1-\delta$, this gives 
$ R^{(m)}(N,\frac bq, f) \gg q^{s}$ with 
$$
s = \eta(m-1)+m\delta-m = m(\eta-(1-\delta)) - \eta
$$
which is positive if $m>\eta/(\eta-(1-\delta))>0$. Thus for $m$
sufficiently large, $ R^{(m)}(N, b/q, f)$ will diverge in these ranges.
\end{proof} 

\section{Proof of Theorem~\ref{thm:4}}\label{sec:6}

Fix $m\ge 2$, some  small $\delta >0$, and let $N$, $q$ be large such that 
$q^{1/3} \leq N \leq q^{m/(m+2) - \delta}$. 
Let $f\geq 0$ is an admissible test
function, $f(0)\neq 0$, and $f_+\geq f$ the smooth majorant appearing in
Lemma~\ref{comparison:lem}. 
We want to show that there exists $b< q$ coprime to $q$ such that 
$R^{(m)}(N,b/q,f_+)$ is large. 

We first  produce $q'$ which is a square, $q'=v^2$, coprime to $q$,
such that 
\begin{equation}\label{size of q'}
q q'\asymp N^3(\log N)^3 \;.
\end{equation} 
To do so, find $v$ in  the interval 
$$
J=\left[ \sqrt{ \frac{N^{3}\log^3 N}{q}},2\sqrt{ \frac{N^{3}\log^3 N}{q}}
\right ]
$$ 
which is coprime to $q$. 
Note that $N^{3}\log^3 N/q \gg (\log q)^3$ since $N\geq q^{1/3}$ 
 and so the existence of such numbers $v$ is assured for any $q$
sufficiently large. Indeed, if $q$ is sufficiently large then 
in any interval $[x,2x]$ with $x\gg (\log q)^{3/2}$ there is a prime
$\ell$ not dividing $q$, since otherwise $q$ would be divisible by all
primes in the interval and consequently $\log q$ would be at least as
large as $\sum_{x\leq p \leq 2x} \log p \sim x$ which contradicts
$x\gg (\log q)^{3/2}$.

We now put $q'=v^2$. Thus $(q',q)=1$ and \eqref{size of q'} holds. 
Because $q'=v^2$ is a square, we may use the divergence principle 
\eqref{blowup:eq} to see that for all $b'$ we have 
$$
R^{(m)}(N,\frac {b'}{q'}, f) \gg \frac 1N (\frac Nv)^m 
$$
Since $v=\sqrt{q'} \ll \sqrt{N^3\log^3 N/q}$, we find that for all $b'$
$$
R^{(m)}(N,\frac {b'}{q'}, f) \gg 
\frac { N^{m-1}q^{m/2} } {N^{3m/2}(\log N)^{3m/2} }
= \frac{q^{m/2}}{N^{m/2+1}(\log N)^{3m/2}} \;.
$$
Now use $q^{1/3}\leq N \leq q^{m/(m+2)-\delta}$  to find that for some
$C>0$, 
\begin{equation}\label{eq:Large'}
R^{(m)}(N,\frac {b'}{q'}, f) \geq C q^{\delta(m/2+1)}(\log q)^{-3m/2}
\end{equation}
uniformly in $b'$ if $q>q_0$. Since $\delta>0$, this diverges with $q$. 

Because $q$, $q'$ are coprime, there are $0<b<q$, $0<b'<q'$ so that 
$bq'-b'q = 1$ and so 
$$
|\frac bq- \frac{b'}{q'}| =\frac 1{qq'} \asymp \frac 1{N^3\log^3 N}   \;.
$$

By the comparison principle (lemma~\ref{comparison:lem}), the two
sequences $\CN = \{\{bn^2/q\}:n\leq N\}$ and  
$\CN' = \{\{b'n^2/q'\}:n\leq N\}$ satisfy 
\begin{equation}\label{eq:comp2}
|R^{(m)}(N,\frac {b}{q}, f)- R^{(m)}(N,\frac {b'}{q'}, f)| \leq 
\epsilon(\CN,\CN') 
R^{(m)}(N,\frac {b}{q}, f_+)
\end{equation}
where $f_+$ is a majorant for $f$, and in particular nonvanishing
at the origin. 
Moreover 
\begin{equation}\label{eq:epsilon}
\begin{split}
\epsilon(\CN,\CN') &= 
N\max_{n\leq N}|\{\frac {bn^2}q\}-\{\frac{b'n^2}{q'}\}| \leq 
|\frac bq- \frac{b'}{q'}| {N^3} \\
&\leq \frac 1{(\log N)^3}
\ll \frac 1{ (\log q)^3}
\end{split}
\end{equation}

We claim that 
$$
R^{(m)}(N,\frac {b'}{q'}, f_+) \geq 
\frac C3 q^{\delta(m/2+1)}(\log q)^{-3m/2}
$$
Indeed, assuming otherwise we have from \eqref{eq:comp2} and
\eqref{eq:epsilon} that 
$$
|R^{(m)}(N,\frac {b}{q}, f)- R^{(m)}(N,\frac {b'}{q'}, f)| = 
o(q^{\delta(m/2+1)}(\log q)^{-3m/2})
$$
which together with \eqref{eq:Large'} forces 
$R^{(m)}(N,{b}/{q},f)>\frac C3 q^{\delta(m/2+1)}(\log
q)^{-3m/2}$. However, since $f_+\geq f\geq 0$ we find that 
$$
R^{(m)}(N,\frac {b}{q},f_+) \geq R^{(m)}(N,\frac {b}{q},f) > 
\frac C3 q^{\delta(m/2+1)}(\log q)^{-3m/2}
$$ 
contradicting our assumption. \qed





\section{Preliminaries on continued fractions} 

We recall the standard notions of the theory of continued fractions
(see e.g. \cite{HW}). 

Given integers $a_0 \in \Z$, $a_1,a_2,\dots \geq 1$, one defines integers 
$p_m$, $q_m$  by the recursion ($m\geq 1$):
\begin{equation*}
\begin{split}
p_m & =a_m p_{m-1} + p_{m-2} \\
q_m & =a_m q_{m-1} + q_{m-2}
\end{split}
\end{equation*}
with $p_{-1} = 1$,  $p_0=a_0$, $q_{-1} = 0$, $q_0=1$.
These satisfy the relations
$$
p_{m}q_{m-1}-p_{m-1}q_{m} = (-1)^{m-1}
$$
and 
$$
p_mq_{m-2}- p_{m-2}q_m = (-1)^m a_m    \;.
$$
The finite continued fraction 
$$
[a_0;a_1,\dots ,a_m]:= a_0 +\cfrac{1}{ a_1+\cfrac{1}{a_2+\cfrac{1}
{\ddots+\cfrac{1}{ a_m }}}}
$$
is then $p_m/q_m$. 

The infinite simple continued fraction $[a_0;a_1,a_2,\dots ]$ is the
limit of the  ``convergents'' ${p_m}/{q_m}$. 
Every irrational $\alpha$ has a  unique continued fraction expansion. 

The convergents give very good rational approximations to
$\alpha$: We have 
$$
\frac 12 \frac 1{q_m q_{m+1} } < |\alpha -\frac {p_m}{q_m}|  < 
\frac 1{q_m q_{m+1}} \;.
$$

The  convergents $p_m/q_m$ are  the ``best'' rational
approximations to $\alpha$, in the following senses: 
If $p/q$ satisfies $|\alpha -p/q|< 1/2q^2$ then $p/q = p_m/q_m$ for
some $m$. Moreover, for $m>1$, if $0<q\leq q_m$ and $p/q \neq p_m/q_m$ then 
$|\alpha- p/q|> |\alpha -p_m/q_m|$.

\section{Proof of Theorem~\ref{thm:2}(a)} 

\subsection{Constructing $\alpha$} 
We want  to find an irrational $\alpha$ such that 
\begin{equation}\label{big R_5}
R^{(5)}(\alpha,N)  \gg \frac{ N^{1/4}}{(\log N)^5} \;.
\end{equation} 
The construction below is due to the referee, who
strengthened and considerably simplified our original argument. 
 
We  construct $\alpha $ by means of its continued fraction
expansion, by inductively finding $a_0, a_1,\dots ,a_m$ so that the
denominators $q_m$ of the  convergents are {\em squares}: $q_m
= v_m^2$. 

To do so, define pairs of integers $(r_m, v_m)$ by $r_{-1} = v_{-1}=0$, $r_0
= v_0=1$, $r_1=v_1=1$ and for $m\geq 1$ 
\begin{equation}\label{recursion of v_m}
v_{m+1} = r_m v_m^2 + v_{m-1}, \qquad r_{m+1} = [\log v_{m+1}] 
\end{equation}
Now set $a_0 = 1$, and for $m\geq 0$ 
$$
a_{m+1} = r_m^2 v_m^2 + 2 r_m v_{m-1} 
$$
Let $\alpha  = [a_0;a_1,a_2,\dots]=[1;1,3,6,\dots]$. 

We claim that the denominator $ q_m$ of  convergent  to
$\alpha$  equals   $v_m^2$. To see 
this, use induction: By the recursion  for the 
convergents, $q_{m+1} = a_{m+1}q_m +q_{m-1}$ and by induction 
\begin{equation*}
\begin{split}
q_{m+1} &= a_{m+1} v_m^2  + v_{m-1}^2 \\
&=        (r_m^2 v_m^2 + 2 r_m v_{m-1})  v_m^2  + v_{m-1}^2 \\
& =       (r_m v_m^2 +v_{m-1})^2 = v_{m+1}^2
\end{split}
\end{equation*}
as required. 

Note also that from the recursion \eqref{recursion of v_m}, 
$$
q_{m+1} \sim r_m^2 q_m^2 \sim q_m^2 (\log q_m)^2
$$
Thus $\alpha$ is of type $ 3+\epsilon$, for all $\epsilon>0$.

Now we want to show that $R^{(5)}(\alpha,N) \gg N^{1/4}/(\log N)^5$. 
Pick $m$ so that  $q_m\leq N <q_{m+1}$. 
We will replace the sequence  of fractional parts $\CN = \{\{\alpha n^2\}:
n\leq N\}$  by a different sequence depending on the size of $N$
relative to $q_m$. 

\subsection{Case 1: Assume that $q_m^{4/3}\leq N < q_{m+1}$} 
Recall that $q_{m+1} \gg q_m^2 \log q_m$ and so this range is 
nonempty. 
Replace $\CN$ by the sequence $\CN' = \{x_n': n\leq N\}$ where $x_n' =
\{p_{m+1} n^2/q_{m+1}\}$. These two sequences have asymptotically equal
correlations since  
$$
|x_n - x_n'| \leq |\alpha -\frac {p_{m+1}}{q_{m+1}}| n^2 < 
\frac {N^2}{q_{m+1} q_{m+2}} \ll \frac {N^2}{r_{m+1}^2 q_{m+1}^3} \ll
\frac 1{N(\log N)^2}
$$
since $q_{m+1}>N$ and $r_{m+1} = [\log q_{m+1}] \gg \log N$. 
Thus by the comparison principle (lemma~\ref{comparison:lem}), it
suffices to work with the new sequence $\CN'$. 

By the divergence principle  (see \eqref{blowup:eq}), since
$q_{m+1}=v_{m+1}^2$,  if the test function $f$ is nonvanishing at the
origin we find that 
$$
R^{(5)}(\CN',f) \gg \frac 1N (\frac N{v_{m+1}})^5 = \frac
{ N^4}{q_{m+1}^{5/2}} \gg \frac{N^4}{r_m^5 q_m^5}   \;.
$$
Since $q_m \ll N^{3/4}$ and $r_m\sim \log N$ we find that 
$$
R^{(5)}(\CN',f) \gg N^{1/4} (\log N)^{-5}
$$
proving \eqref{big R_5} when $q_m^{4/3}< N< q_{m+1}$. 

\subsection{Case 2: $q_m\leq N< q_m^{4/3}$}   
Set 
$$ 
M=\frac {q_m^{3/2}}{N^{1/2}}
$$  
Note that since $q_m \leq N<q_m^{4/3}$, $M$ lies between $N^{5/8}$ and
$N$. 
We replace $\CN$ by the sequence $\CN''=\{x_n'': n\leq N\}$ where 
$$
x_n'' = \begin{cases} \{\frac {p_m}{q_m} n^2\}, & n\leq M \\
                     x_n = \{\alpha n^2\}, & M<n\leq N
\end{cases}
$$

To check that correlations of $\CN$ and $\CN''$ are asymptotically
equal, we need to see that $|x_n-x_n''|  = o(1/N)$. For $n>M$ this
certainly holds, while for $n\leq M$ we have 
$$
|x_n-x_n''| \leq |\alpha -\frac {p_m}{q_m}| n^2 \leq 
\frac {M^2}{q_m q_{m+1}}
$$
Now use $q_{m+1} \gg r_m^2 q_m^2$ and since $q_m \sim (M^2 N)^{1/3}$
and $r_m \gg \log N$, 
we have 
$$
q_m q_{m+1} > r_m^2 q_m^3 \gg (\log N)^2 M^2 N
$$
which gives $|x_n-x_n''| \ll 1/N(\log N)^2$ as required.  

Now we study the sequence $\CN''$. The number $0$ occurs in $\CN''$ 
if $n\leq M$ is a multiple of $v_m$: $n=v_m n'$, since then 
$x_n'' = \{\frac{p_m}{q_m}n^2 \} = \{p_m (n')^2\} = 0$. Thus $\vec 0$
occurs as a difference of $5$-tuples of elements of $\CN''$ at least
$\gg [M/v_m]^{5}$ times. Thus if the origin lies in the support of the
test function $f$ then 
$$
R^{(5)}(\CN'',f) \gg \frac 1N (\frac M{v_m})^5  = \frac
{M^5}{Nq_m^{5/2}} 
$$
Since $M^2 = q_m^3/N$, and $ N<q_m^{4/3}$ we get
$$
R^{(5)}(\CN'',f) \gg  \frac{q_m^5}{N^{7/2}} \gg N^{15/4 - 7/2} = N^{1/4}
$$

\section{Proof of Theorem~\ref{thm:2}(b)} 
Let $\sigma >23/8$. 
We construct $\alpha = [a_0;a_1,\dots,a_m,\dots]$ 
which will be of  type
$\sigma$   and for which $\limsup R^{(5)}(\alpha, N) = \infty$ 
by an inductive construction of the partial quotients $a_m$. 

Suppose we have already found $a_0,\dots ,a_{m-1}$, from which we got
the partial convergents $p_j/q_j$, $j=0, \dots ,m-1$. 
Now take an integer $\ell \sim q_{m-1}^{(\sigma-2)/2}$, which is
coprime to $q_{m-1}$ (this is certainly possible for $m\gg 1$, say
take $\ell$ a prime between $ q_{m-1}^{(\sigma-2)/2}$ and 
$2 q_{m-1}^{(\sigma-2)/2}$ which does not divide $q_{m-1}$). 
Also set $v_m = \ell$. 

Because $\ell$ and $q_{m-1}$ are coprime,  there is a  unique solution
$t=a_m$ of the congruence 
\begin{equation}\label{congruence mod l^2}
tq_{m-1}+q_{m-2} = 0 \mod \ell^2
\end{equation}
which lies in $[\ell^2, 2\ell^2)$. 
Then $a_m \sim \ell^2 \sim q_m^{\sigma-2}$, 
and 
$$
q_{m} :=a_m q_{m-1}+q_{m-2} \sim q_{m-1}^{\sigma-1}     \;.
$$ 
Thus $\alpha$ is of type  $\sigma+\epsilon$ for all 
$\epsilon>0$. 
Moreover $q_m$ is divisible by $v_m^2$ by \eqref{congruence mod l^2}. Thus 
$$
q_m = u_m v_m^2
$$
for some integer $u_m$, and 
$$
v_m\sim q_{m-1}^{(\sigma-2)/2} \sim q_m^{(\sigma-2)/(2\sigma -2)} 
$$

Now take 
$$
N_m \sim \frac {q_m^{\sigma/3} }{ \log q_m}
$$
We will see that $R^{(5)}(\alpha, N_m)\to \infty $ as $m\to \infty$. 

To see this, note that in the sequence 
$\{\alpha n^2: n\leq N_m \}$ we may replace 
 $\alpha$ by the partial convergent $p_m/q_m$ without changing the
limiting correlations. To see this, note that by
Lemma~\ref{comparison:lem} it
suffices to check that $1/q_mq_{m+1} = o(1/N_m^3)$. Indeed, we have 
$$
\frac 1{q_m q_{m+1}} \sim \frac 1{q_m^\sigma} \ll \frac 1{N_m^3 (\log
N_m)^3}
$$
as required. 

To see that $R^{(5)}(N_m, \frac {p_m}{q_m}, f)$  diverges for positive
test functions $f$ with $f(0)\neq 0$, 
use the divergence principle \eqref{blowup:eq} to find 
$$
R^{(5)}(N_m, \frac {p_m}{q_m}, f) 
\gg \frac 1{N_m} (\frac {N_m v_m}{q_m})^5 \;.
$$
Now use $N_m \sim q_m^{\sigma/3}/\log q_m$ and 
$v_m\sim q_m^{(\sigma-2)/(2\sigma -2)}$ to find 
$$
R^{(5)}(N_m, \frac {p_m}{q_m}, f) \gg \frac {q_m^{E}}{(\log q_m)^5}
$$
where 
$$
E = \frac {4\sigma}3 + 5\frac{\sigma - 2}{2\sigma -2} - 5 = 
\frac {\sigma(8\sigma-23)}{6(\sigma -1)} \;. 
$$
Since $\sigma>23/8$, we have $E>0$ which gives divergence of $R^{(5)}$. \qed

\section{Proof of Theorem~\ref{thm:3}}\label{sec:3}

Fix $m\ge 2$, $f$ and $\delta >0$. By approximating $f(0,x_2,\cdots
,x_m)$ from above and below with step functions we see that it is
enough to prove the statement  for a function $f$ symmetric, satisfying 
$f(x+(t,t,\cdots ,t))=f(x)$ for all $t\in {\bf R}$ and such that 
$f(0,x_2,\cdots ,x_m)$ is the characteristic function of a nice compact
set $I\subset {\bf R }^{m-1}$. In other words, given such an $I$ and 
$m,\delta $ as above, it is enough to show that as $q\rightarrow \infty $
one has 
\begin{equation}\label{e31}
R^{(m)}(N,b/q,I)\rightarrow Vol(I)
\end{equation}
uniformly for $(b,q)=1$ and  $q^{1-\frac {1}{2m}+\delta }\le N\le 
\frac {q}{\log q}$, where $NR^{(m)}(N,b/q,I)$ is the number of 
tuples $(x_1,\cdots ,x_m)$ with distinct components $x_1,\cdots ,x_m$
in $ \{1,\cdots ,N\}$ such that 

$$
N(\{\frac {bx_1^2}{q}\}-\{\frac {bx_2^2}{q}\},\cdots ,
\{\frac {bx_{m-1}^2}{q}\}-\{\frac {bx_m^2}{q}\})\in I.
$$
Given a large prime number $q$ and $b,N$ as above, we write $R^{(m)}(N,b/q,I)$
in the form 
$$
R^{(m)}(N,b/q,I)=\frac 1N\sum _{\vec a \in sI}^* \nu (N,\vec a)
$$ 
where $s=\frac qN $ is the dilate factor, and 
$$
\nu(N,\vec a)=\# \{ 1\le x_i \le N: bx_i^2-bx_{i+1}^2=a_i\pmod q,
1\le i\le m-1 \} \;.
$$
Here $\sum ^*$ means the summation is over the vectors $\vec a $
for which the partial sums $A_i=\sum _{k\ge i}a_k,A_m=0$,
are distinct, a condition which comes  from the 
requirement that the $m$-tuples $x=(x_1,... ,x_m)$ to be counted
in $R^{(m)}(N,b/q,I)$ have distinct components.
Let 
$$
h_{\vec a}(\vec x) = 
\begin{cases} 1, & b(x_i^2-x_{i+1}^2)=a_j \pmod q  , \qquad i=1,\cdots ,
m-1\\
0 & \text{else.}
\end{cases}
$$
Thus:
$$
\nu (N,\vec a)= \sum _{1\le x_1,\cdots ,x_m\le N} h_{\vec a} (\vec x) \;.
$$
We now use the Fourier expansion:
$$
\nu(N,\vec a)=\sum _{\vec r \pmod q}\hat h_{\vec a}(\vec r) \prod 
_{i=1}^m F_N(r_i)
$$
where
$$
\hat h_{\vec a}(\vec r)=\frac {1}{q^m} \sum _{\vec y \pmod q}h_{\vec a}
(\vec y)e\Big (-\frac {{\vec r}\cdot {\vec y}}{q}\Big )
$$
and:
$$
F_N(r_i)=\sum_{1\le x_i\le N} e\Big (\frac {r_ix_i}{q}\Big ).
$$
These last sums are geometric series which can be bounded by:
\begin{equation}\label{e32}
\|F_N(r_i)\| \ll \min \{ N,\frac {q}{|r_i|} \}
\end{equation}

where the residues $r_i$ are assumed to lie in the interval
$[\frac {-q}{2},\frac q2]$.
In 
$$
R^{(m)}(N,b/q,I)=\frac 1N \sum _{\vec a \in sI} \sum _{\vec r \pmod q }^*
\hat h_{\vec a} (\vec r) \prod _{i=1}^m F_N(r_i)
$$
we isolate the contribution of $\vec r=0$ to get the main term :

\begin{equation}\label{e33}
R ^{(m)}(N,b/q,I)=\mathcal M+ \mathcal E
\end{equation}
with

\begin{equation}\label{e34}
\mathcal M=N^{m-1} \sum _{\vec a \in sI}^*\hat h_{\vec a}(0)
\end{equation}
and
\begin{equation}\label{e35}
\mathcal E=\frac 1N \sum _{0 \ne \vec r \pmod q } \prod _{i=1}^mF_N(r_i)
\sum _{\vec a \in sI}^* \hat h_a(\vec r).
\end{equation}

We first estimate the main term. 
For any $\vec a$  
let $C(\vec a,q)$ be the curve mod $q$ given by the system of
congruences:
$$
\begin{array}{l}
bx_1^2-bx_{2}^2=a_1 \pmod q \\
\cdots \\
bx_{m-1}^2-bx_{m}^2=a_{m-1} \pmod q. \end{array}
$$
One has $\hat h_{\vec a }(0)=
\frac {1}{q^m}\nu(\vec a ,q)$, where 
$\nu(\vec a ,q)$
is the number of points on the curve  $C(\vec a,q)$. Thus 

\begin{equation}\label{e36}
\mathcal M=\frac {N^{m-1}}{q^m}\sum _{\vec a \in sI}^* \nu (\vec a,q).
\end{equation}
We want to show that as $q\rightarrow \infty $ one has:
\begin{equation}\label{e37}
\mathcal M=Vol(I)+o(1).
\end{equation}
For any $\vec a =(a_1,\cdots ,a_{m-1})$ denote by 
$r_{eff}(\vec a,q)$ the number of distinct $y_j$ satisfying the following 
system:
\begin{equation}\label{e38}
y_i-y_{i+1}=a_i(\mod q),1\le i\le m-1.
\end{equation}
Since the solutions of the homogeneous system
$$
y_i-y_{i+1}=0(\mod q),  1\le i\le m-1
$$
are spanned by $(1,\cdots ,1)$, $r_{eff}(\vec a,q)$ is well-defined 
(independent of the particular solution $y$ of (\ref{e38})).
Using the Riemann Hypothesis for curves over finite fields (Weil \cite{We})
one obtains (see also \cite{K-R}, Proposition 4):

\begin{equation}\label{e39}
\nu(\vec a,q)=2^{m-r_{eff}(a,q)}(q+B(\vec a,q))
\end{equation}
with 
\begin{equation}\label{e310}
|B(\vec a,q)|\ll_mq^{\frac12}.
\end{equation}
We define roots $\sigma _{ij}(\vec a)$,$1\le i<j\le m$ by
\begin{equation}\label{e311}
\sigma _{ij}(\vec a)=\sum _{k=i}^{j-1} a_k
\end{equation}
so that $\sigma _{i,i+1}(\vec a)=a_i$,
 $ \sigma _{ij}=\sum _{k=i}^{j-1} \sigma _{k,k+1}$.
We set $D(\vec a)=\prod _{1\le i\le j\le m}\sigma _{ij}(\vec a)$.
The solutions of (\ref{e38}) are all distinct (i.e. $r_{eff}(\vec a ,q)=m$)
if and only if $q$ does not divide $D(\vec a )$, since 
$y_i-y_j=\sum _{k=i}^{j-1}y_k-y_{k+1}=\sum _{k=i}^{j-1}a_k=
\sigma _{ij}(\vec a )$.
Note that $D(\vec a)$ is a nonzero integer for any $\vec a $ which appears 
 in the above summations $\sum _{a\in sI}^*$.
In our case $q$ does not divide $D(\vec a )$, since for $N$ large enough in terms of
 $I$ each factor $\sigma _{i,j}(\vec a )$ of $D(\vec a)$ is in absolute value smaller
 than $q$. Therefore $r_{eff}(\vec a ,q)=m$ and (\ref{e39}) and (\ref{e310}) give

\begin{equation}\label{e312}
\nu (\vec a,q)=q+O_m(q^{\frac 12})
\end{equation}
for all $\vec a $ which appear in (\ref{e36}). Then (\ref{e36}) implies that

\begin{equation}\label{e313}
\mathcal M =\frac {N^{m-1}}{q^m} 
(q+O_m(q^{\frac 12}))\sum _{a\in sI}^*1=
\frac {1}{s^{m-1}}(1+O_m(\frac {1}{q^{\frac 12}}))\sum _{a\in sI}^*1.
\end{equation}
The number of integer points $\vec a \in sI $ which lie in the union of the 
hyper-planes $\sigma _{ij}(\vec a )=0$ is 
$O_{m,I}(s^{m-2})$, while by the Lipschitz principle (see Davenport \cite{Da})
it follows that:
$$
\#(sI\cap {\bf Z}^{m-1})=s^{m-1}Vol(I)+O_{m,I}(s^{m-2}).
$$
Therefore:
\begin{equation}\label{e314}
\sum _{a\in sI}^*1=\#(sI\cap {\bf Z}^{m-1})-\#\{\vec a \in sI : D(\vec a)=0\}
\end{equation}
$$
=
s^{m-1}Vol(I)+O_{m,I}(s^{m-2}) 
$$
and from (\ref{e313}) we get
$$
\mathcal M=(1+O_m(\frac 1{\sqrt q}))(1+O_{m,I}(\frac 1s))
$$
which proves (\ref{e37}).

We now proceed to estimate the remainder $\mathcal E$.
For any $\vec a$  and $ \vec r $ we have:
$$
\hat h_{\vec a}(\vec r)=\frac {1}{q^m}\sum _{\vec y \in C (\vec a ,q)}
e\Big (-\frac {\vec r \cdot \vec y }{q}\Big ).
$$
Applying Weil's Riemann Hypothesis for curves over finite fields
one has (see \cite{Bo}, Theorem 6)
\begin{equation}\label{e315}
\Big | \sum_{y \in C(a,q)} e\Big ( -\frac {\vec r \cdot \vec y}q \Big )
\Big | \ll_m \sqrt q
\end{equation}
unless the linear form $\vec r \cdot \vec y$ is constant along the curve.
For $\vec a$ as in (\ref{e35}) this only happens if $\vec r=0$. For,
let  $\vec r \neq 0$ be such that $\vec r \cdot \vec y$ 
is constant along the curve. Then, in the function field
$\bar {\bf F}_q(Y_1,\cdots,Y_m)$ of the curve, where $\bar {\bf F}_q$
denotes the algebraic closure of ${\bf F}_q= {\bf Z}/q{\bf Z}$, $Y_1$
is a variable and $Y_2,\cdots ,Y_m$ are algebraic functions such that
$$
Y_i^2 = Y_1^2 -\frac {a_1+\cdots +a_{i-1}}b
$$
for $2 \le i \le m$, we will have an equality  $\vec r \cdot \vec Y = c$,
with $c \in  \bar {\bf F}_q$. If we choose $j_0 \in \{1,\cdots,m\}$ such that
$r_{j_0} \neq 0$ then $Y_{j_0}$ will lie in 
$\bar {\bf F}_q(Y_1,\cdots,Y_{j_0-1},Y_{j_0+1},\cdots, Y_m)$
and hence
\begin{equation}\label{e316}
\bar {\bf F}_q(Y_1,\cdots,Y_m)=
\bar {\bf F}_q(Y_1,\cdots,Y_{j_0-1},Y_{j_0+1},\cdots, Y_m).
\end{equation}
Now for any unique factorization domain $D$ of characteristic $\neq 2$
and any distinct primes $p_1,\cdots,p_t$ in $D$ one has
$$
[K(\sqrt {p_1},\cdots,\sqrt {p_t}):K]=2^t
$$
where $K$ denotes the quotient field of $D$ (see Besicovitch 
\cite{Be})).
Applying this with $D=\bar {\bf F}_q[Y_1]$ and  
$p_i = Y_1^2-\frac {a_1+\cdots +a_{i}}b $ for $1 \le i \le m-1$ we get:
$$
[\bar {\bf F}_q(Y_1,\cdots,Y_m):\bar {\bf F}_q(Y_1)]=2^{m-1}.
$$
By the same argument we see that 
$$
[\bar {\bf F}_q(Y_1,\cdots,Y_{j_0-1},Y_{j_0+1},\cdots, Y_m)
:\bar {\bf F}_q(Y_1)]=2^{m-2}
$$
which contradicts (\ref{e316}). It follows that for all $\vec r$ and 
$\vec a$ which appear in (\ref{e35}), the inequality (\ref{e315}) holds true
and one has:
$$
|\hat h_{\vec a}(\vec r)| \ll_m \frac 1{q^{m-\frac 12}} .
$$
This implies that
\begin{equation}\label{e317}
|\mathcal E| \ll_m \frac 1{Nq^{m-\frac 12}} \sum_{0 \neq \vec r \pmod q}
\Big ( \prod_{i=1}^m |F_N(r_i)|\Big ) \sum_{\vec a \in sI}^*1 .
\end{equation}
We use (\ref{e32}) and (\ref{e314}) in (\ref{e317}) to conclude that
\begin{equation}\label{e318}
|\mathcal E| \ll_{m,I} \frac {s^{m-1}}{Nq^{m-\frac 12}} 
\sum_{\vec r \pmod q}
\prod_{i=1}^m \min \{N,\frac q{|r_i|}\}
\end{equation}
$$
\ll_m \frac {q^{m-\frac 12}\log ^m q}{N^m} \le 
\Big (\frac {\log q}{q^{\delta}}\Big )^m.
$$
The theorem now follows from 
(\ref{e33}), (\ref{e37}) and (\ref{e318}).

\appendix
\section{Square factors of rational approximants}\label{app:A}

Let $\alpha$ be a real number and $a_n/q_n$ a sequence of rational
approximants of $\alpha$: $|\alpha-a_n/q_n|<1/q_n^2$, and
$q_n\to\infty$.
In view of Theorem 1',
we want to investigate the square parts of the denominators $q_n$,
keeping in mind that large square parts rule out Poisson statistics
for the correlation functions.

\begin{defn}
A sequence $\{q_n\}$ is {\em almost square-free} if  $\forall
\epsilon>0$, all square divisors $s_n^2$ of $q_n$ satisfy
$s_n\ll_\epsilon q_n^\epsilon$.
\end{defn}

\subsection{A metric result}
We will show that for almost all $\alpha$, we have: 
If $a_n/q_n$ is a sequence of
rational approximants of $\alpha$ (that is $|\alpha-a_n/q_n|<1/q_n^2$, and
$q_n\to\infty$), then $\{q_n\}$ is almost square-free.

In fact, we show  more: 
For an integer $q\geq 1$, we write $q=\tq s^2$ with $\tq$ square-free.
Let $\SF$ be the set of  integers $q$ whose largest  square factor
$s^2$ satisfies  $s\leq \log^2 \tq$. 
We will show that almost all reals $\alpha$ have rational approximants
whose denominators are in $\SF$ except for finitely many exceptions.

\begin{prop}\label{prop:sf}
For all  reals $\alpha$ outside  a set of measure zero, there is a
$Q=Q(\alpha)>1$ so that if
$|\alpha-a/q|<1/q^2$ and $q\geq Q$  then $q\in \SF$. 
\end{prop}
The proof of this  follows from a well-known
general principle: Given a sequence  of integers $\CN$, we say that a
real number $\alpha$ is $\CN$-approximable if there are {\em 
infinitely many}
rationals $a/q\neq \alpha$ with denominator $q\in\CN$ and
$|\alpha-a/q|<1/q^2$. For instance, we may take as $\CN$ the complement
of $\SF$. To prove Proposition~\ref{prop:sf}, we will use
\begin{lem} \label{genprin}
Suppose that $\CN$ is a sequence such that
$$\sum_{q\in \CN}\frac 1q <\infty .$$
Then the set of $\CN$-approximable reals has measure zero.
\end{lem}
\begin{proof} Without loss of generality we
will assume that $0<\alpha<1$. For each pair of coprime integers
$(a,q)$ with $1\leq a<q$, denote by $I_{a,q}$ the interval
$$I_{a,q} = (\frac aq-\frac 1{q^2},\frac aq+\frac 1{q^2} )$$
Then $\alpha$ is $\CN$-approximable if and only if  it lies in infinitely many of
the intervals $I_{a,q}$ with $q\in \CN$. That is for all $N\geq 1$,
$\alpha$ lies in
$$M_N:=\cup_{N\leq q\in \CN}\cup_{1\leq a<q}I_{a,q}.$$
Thus we need to compute the measure of $M:=\cap_{N\geq 1} M_N$. Since
$M_{N}\supseteq M_{N+1}\supseteq\dots$, we have
$$
\meas(M) = \lim_N \meas(M_N) \leq \lim_N \sum_{N\leq q\in \CN} \sum_{a=1}^q
\meas(I_{a,q}) \ll \lim_N \sum_{N\leq q\in \CN} \frac 1q
$$
(allowing overlap of the intervals).  Since
$\sum_{q\in\CN}1/q<\infty$, the above limit is zero. 
\end{proof}

Thus to prove Proposition~\ref{prop:sf}, it suffices to show  
\begin{equation*}
\sum_{q\notin \SF} \frac 1q <\infty \;.
\end{equation*}
We rewrite this sum by grouping together those $q$ with the same
square-free kernel $\tq$: Writing $q=fm^2$, $\tq=f$, then
\begin{equation*}
\sum_{q\not\in \SF} \frac 1q = \sum_{f\mbox{ square-free}}
\sum_{\substack{\tq=f \\ q\not\in \SF}} \frac 1q
=\sum_{f\mbox{ square-free}} \frac 1f
\sum_{q=fm^2 \not\in \SF}\frac 1{m^2}
\end{equation*}
Now if $q\not\in \SF$, $\tq=f$ then $m>\log^2 f$. Thus for each $f$,
$$
\sum_{q=fm^2 \not\in \SF}\frac 1{m^2} =\sum_{m>\log^2 f}\frac
1{m^2} \ll \frac 1{\log^2 f}
$$
and so
$$
\sum_{q\not\in \SF} \frac 1q  \ll \sum_{f\mbox{ square-free}}
\frac 1f \frac 1{\log^2 f }<\infty
$$
as  required. \qed

\subsection{Algebraic $\alpha$}
For real {\em algebraic} $\alpha$, the analogue of 
Proposition~\ref{prop:sf} 
follows from a standard belief in diophantine analysis, namely the
``ABC Conjecture'' of Masser and Oesterle: 
Define the {\em radical} of an integer $N$ as the
product of all primes dividing it: $\rad(N):= \prod_{p\mid N} p$.
The ABC conjecture is the assertion that whenever we have an equation
in coprime integers $A+B+C=0$, then
\begin{equation}\label{abc conj} 
|A| \ll_\epsilon \rad(ABC)^{1+\epsilon}
\end{equation}
for all $\epsilon>0$.
This implies a seemingly stronger statement:
Suppose that $G(x,y)\in \Z[x,y]$ is a homogeneous form with integer
coefficients and no repeated factors, and $m,n$ coprime integers. 
Then for all $\epsilon>0$
\begin{equation}\label{super ABC}
\max(|m|,|n|)\sp {\deg(G)-2-\varepsilon}\ll_\epsilon \rad(G(m,n))\;,
\end{equation}
where $\deg(G)$ is the degree of $G$.
The deduction of \eqref{super ABC} from \eqref{abc conj} and 
a theorem of  Belyi \cite{Belyi} was noted by Elkies \cite{Elkies}
and by  Langevin \cite{Langevin}.  
The ABC-conjecture \eqref{abc conj} is the
special case of the ternary form $G(x,y)=xy(x+y)$.

The corollary \eqref{super ABC} of the ABC conjecture implies the
analogue of Proposition~\ref{prop:sf} for irrational algebraic $\alpha$.
Indeed, let $f(x)$ be the minimal polynomial of $\alpha$, of degree
$d>1$, and write $f(x/y) =  F(x,y)/y^d$ with $F(x,y)\in \Z[x,y]$.
Suppose that $p/q$ is an approximant of $\alpha$:
$|\alpha -p/q|<1/q^2$, with $p$, $q$ coprime. 
Since $f(x)$ is irreducible, $f'(\alpha)\neq 0$
and thus by the mean value theorem, for some $\xi$ between $\alpha$
and $p/q$,
$$
|f(\frac pq)|  =|f(\frac pq) - f(\alpha)|  = |\alpha - \frac
pq||f'(\xi)| \ll \frac 1{q^2}
$$
On the other hand,
$$ f(\frac pq) = \frac{F(p,q)}{q^d}$$
and so we find
$$
|F(p,q)|\ll q^{d-2} \;.
$$

By \eqref{super ABC}, taking $G(x,y) = xy F(x,y)$ and noting that
$|p|\ll q$, we get for all $\epsilon>0$
$$
q^{d-\epsilon} \ll_\epsilon \rad(pqF(p,q))\leq |pF(p,q)|\rad(q) \ll
q^{d-1}\rad(q) \;.
$$
Thus if $q=\tilde q s^2$ then
$$
(\tilde q s^2)^{d-\epsilon} \ll_\epsilon \rad(\tilde q s) \leq \tilde
q s
$$
and so $s\ll_\epsilon q^\epsilon$.

\end{document}